\documentclass[12pt]{article}
\usepackage[T2A]{fontenc}
\usepackage[cp1251]{inputenc}
\usepackage[russian]{babel}
\usepackage{amssymb,amsmath,latexsym,amsthm}
\usepackage{indentfirst}
\usepackage{verbatim}
\theoremstyle{plain}
\theoremstyle{remark}

\title{ZERO (SUB-)SEQUENCES OF ENTIRE FUNCTIONS}
\author{\bf B.\,N.\,Khabibullin, G.\,R.\,Talipova, F.\,B.\,Khabibullin}
\date{{Bashkir State University,
Khabib-Bulat@mail.ru}}
\begin{document}
\maketitle  

\noindent
{{\bf 1.} \sc Zero subsets for wieghted classes.} Let $D$ be a domain of the complex plane $\mathbb C$.  Denote by $\text{sbh\,}(D)$ the class of subharmonic functions in $D$, and $\text{sbh\,}(D) \ni \boldsymbol{-\infty}\colon z\mapsto -\infty$, $z\in D$ {\large(}see [1]{\large)}.

A function $V\in \text{sbh\,}\bigl(\mathbb C \setminus \{0\}\bigr)$  is called a Jensen potential if 
 there exists a $R_V>0$ such that $0\leq V(z) \leq \max\bigl\{0, \log \bigl({R_V}/{|z|}\bigr)\bigr\}$ for all $z\in \mathbb C\setminus \{0\}$.
The class of all Jensen potentials will be denoted by $PJ_0$.
Let $M\in \text{sbh\,}(\mathbb C )$ with $M\neq \boldsymbol{-\infty}$, and let $\nu_M := \frac{1}{2\pi}\Delta M$ be the Riesz measure for M, where 
$\Delta$ is the Laplace operator on distributions. Let $ {\tt Z} =\{{\tt z}_k\}_{k=1,2,}\subset \mathbb C\setminus \{0\}$ be a  sequence of points.

\medskip
\noindent
\textbf{Theorem 1.} {\it  If there are a non-zero entire function\/ $f$ {\rm {\large(}we write $f\neq 0${\large)}} vanishes on\/ $\tt Z$ {\rm {\large(}we write 
$f({\tt Z})=0${\large)}} and $r_0>0$ such that\/ $|f|\leq e^M$ on\/ $\mathbb C$ (pointwise), and
$\inf_{|z|<r_0} M(z)>-\infty$,  then 
$$
\sup \left\{\sum_{k} V({\tt z}_k)-\int_{\mathbb C \setminus \{0\}} V \,d\nu_M \colon  V\in PJ_0 \right\}<+\infty.
$$
Conversely, if 
$$
\sup \left\{\sum_{k} V({\tt z}_k)-\int_{\mathbb C \setminus \{0\}} V \,d\nu_M \colon  V\in PJ_0 \cap C^{\infty}\bigl( \mathbb C \setminus \{0\}\bigr)\right\}<+\infty,
$$
then, for each number $N> 0$, there is an  entire function\/ $f\neq 0$ such that $f({\tt Z})=0$ and
	$$
	\bigl|f(z)\bigr|\leq \exp \left(\frac{1}{2\pi}\int_0^{2\pi} M\Bigl(z+
	\frac{1}{1+|z|^{N}}\,e^{i\theta}\Bigr) \, d \theta \right) \quad \text{for all $z\in \mathbb C$}. 
	$$}
	\medskip
	
 This Theorem 1 improves our earlier results from [2]--[8] etc. 

\noindent
{{\bf 2.} \sc Subsequences of zeros for  functions with majorant of Cartwright class.}
A function  $M\in \text{sbh\,}(\mathbb C)$  belong to the (Cartright) class $\mathcal C$ {\large(}see [9]{\large)}
if  $M$ is harmonic outside the real axis $\mathbb R$,  $M(0)=0$, $M(z)=M(\bar z)$ for all $z\in \mathbb C$, and
$$
\limsup_{z\to \infty}\frac{M(z)}{|z|}<+\infty, \quad \int_{\mathbb R}\frac{\max\bigl\{0, M(x)\bigr\}}{x^{2}}\,d x<+\infty.
$$ 
Let $\nu_M$ be the Riesz measure of $M\in \mathcal C$. We denote by  $\nu_M^{\mathbb R}$ 
 its distribution function 
\begin{equation*}\label{nuR}
\nu_M^{\mathbb R}(t):=\begin{cases}	
-\nu_M\bigl( [t,0) \bigr) \;  &\text{ when } t<0 ,\\
\nu_M\bigl( [0,t] \bigr)  \;  &\text{ when } t\geq 0,
\end{cases}
\qquad t\in \mathbb R.
\end{equation*}
We introduce the class $R\mathcal P_0$ of test functions as the subclass of all upper semicontinuous functions $\phi \colon \mathbb R \setminus \{0\} \to [0,+\infty)$ such that  
$$
\text{$ \phi (x)\equiv 0$ for $ |x|\geq R_\phi$} , \quad   \limsup\limits_{0\neq x\to 0} \dfrac{\phi (x)}{-\log |x|}\leq 1,
$$
and, for each   $x_0\in \mathbb R\setminus \{0\}$, there is $r_0\in \bigl(0, |x_0|\bigr)$ such that 
\begin{equation*}
\phi (x_0)\leq \frac{1}{\pi^2}  \int^{+\infty}_{-\infty} \phi(x_0+x)\,\frac{1}{x}\, \log \left|\frac{x+r}{x-r}\right| dx \quad \text{for all  $r\in \bigl(0, r_0\bigr)$}. 
\end{equation*}

\medskip
\noindent
\textbf{Theorem 2.} {\it Let $M\in \mathcal C$, ${\tt Z}=\{{\tt z}_k\}_{k=1,2,\dots}\subset \mathbb C \setminus \{0\}$. 

If there exists an entire function $f\neq 0$ such that $f({\tt Z})=0$ and $|f|\leq e^M$ on $\mathbb C$ (pointwise), then 
\begin{equation*}
			\sup \left\{ \sum_{k} ({\mathrm P} \phi) ({\tt z}_k)
-\int_{-\infty}^{+\infty}\phi(t) \, d \nu_{M}^{\mathbb R}(t) \colon {\phi \in \,R\mathcal P_0}\right\}	<+\infty ,
\end{equation*}
where $(\mathrm P\phi) ({\tt z}):=\phi ({\tt z})$ for $\Im {\tt z}=0$, and, for $\Im {\tt z} \neq 0$,
\begin{equation*}
	(\mathrm P\phi) ({\tt z}):=	\dfrac{1}{\pi}\int_{-\infty}^{+\infty} \Bigl|\Im \dfrac{1}{x-{\tt z}}\Bigr| \phi (x) \, dx 
	\quad \text{\rm (the Poisson integral).}
\end{equation*}

Conversely, if\/ 
\begin{equation*}
	\sup \left\{ \sum_{k} ({\mathrm P} \phi) ({\tt z}_k)
-\int\limits_{-\infty}^{+\infty}\phi(t) \, d \nu_{M}^{\mathbb R}(t) \colon {\phi \in \,R\mathcal P_0}\cap C^{\infty}\bigl(\mathbb R \setminus \{0\}\bigr)\right\}	<+\infty ,
\end{equation*}
then, for each $N>0$, there exists an entire function $f\neq 0$ such that $f({\tt Z})=0$ and 
$|f|\leq \exp M_{N}$ on $\mathbb C$ where
$$
M_{N}(z):=\begin{cases}
\frac{1}{2\pi}\,\int\limits_{0}^{2\pi}M\bigl(z+\frac{1}{1+|\Re z|^N}e^{i\theta}\bigr) \, d \theta\quad
&\text{when $|\Im z|\leq |\Re z|^{-N}$},\\
M(z)\quad  &\text{when $|\Im z|\geq |\Re z|^{-N}$} .
\end{cases}
$$
}

 This Theorem 2 generalizes our previous  results  [10]--[11].

\medskip
\noindent
{{\bf 3.} \sc (Non-)uniqueness sequences.} 
For a subset $S\subset \mathbb C$, we denote by $\text{sbh} (S)$ the class  of functions that  are subharmonic on some open set containing $S$. We set 
$$
\text{sbh}^+(S):=\bigr\{u\in \text{sbh}(S)\colon u(z)\geq 0 \text{ for all  $z\in S$}\bigl\}.
$$
Given $r>0$ and $b>0$, we set 
$${D (r)}:=\bigl\{ z \in \mathbb C\colon  |z| < r\bigr\},\quad 
\overline{D (r)}:=\bigl\{ z \in \mathbb C\colon  |z| \leq r\bigr\},
$$ 
\begin{equation*}
	\text{sbh}_0^+(r;\leq b):= \left\{ v\in \text{sbh}^+\bigl(\mathbb C\setminus {{D(r)}}\,\bigr)\colon 
	\lim_{z\to \infty} v(z){=}0, \; \sup\limits_{|z|=r}v(z)\leq b \right\}.
\end{equation*}
We denote by  $\text{const}_{a_1, a_2, \dots}$  a constant depending only on $a_1,a_2, \dots$.
\medskip

\noindent
\textbf{Theorem 3.} {\it Let  $M\in \text{\rm sbh\,}(\mathbb C)$ with the Riesz measure $\nu_M$,  $r_0,b>0$. Then 
	there are numbers $C:=\text{\rm const}_{r_0,b}>0$,  $\overline{C}_M:=\text{\rm const}_{r_0,   M}\geq 0$
such that, for any  $v{\in}  \text{\rm sbh}_0^+(r_0;\leq b) $ and for each function $u\in \text{\rm sbh\,} (\mathbb C)\setminus \{\boldsymbol{-\infty}\}$
with the Riesz measure $\nu_u$, the pointwise  inequality $u\leq M$ on $\mathbb C\setminus D(r_0)$ entails the inequality\/ 
\begin{equation*}
	\int_{\mathbb C\setminus D(r_0)}  v \,d {\nu}_u 		\leq	\int_{\mathbb C\setminus D(r_0)}  v \,d {\nu}_M	+C\, \overline{C}_M-C u(z_0),
\tag{\rm C}
\end{equation*}
where  a constant $\overline{C}_M<+\infty$ is positively homog\-e\-n\-e\-o\-us of $M$, i.\,e.
 $\overline{C}_{aM}=a\overline{C}_M$ for  $a\in [0,+\infty)$, and   upper semi-additive of $M$, i.\,e.
		$\overline{C}_{M_1+M_2}\leq \overline{C}_{M_1}+\overline{C}_{M_2}$ for  $M_1, M_2\in \text{\rm sbh\,}(\mathbb C)$.
Besides, the integral in the right parts of  the inequality\/ {\rm (C)}   can be replaced by the  integral 		$\int\limits_{\mathbb C\setminus D(r_0)} M\, d \nu_v$, where $\nu_v$ is the Riesz measure of $v$.
}
\medskip 

\noindent
\textsc{Corollary 1.} 
{\it Let  $M\in \text{\rm sbh} (\mathbb C)$ with the Riesz measure $\nu_M$, and 
$ f$ be an entire function. Suppose that  there is a number $r_0>0$ such that 
$\log |f|\leq M$ on $\mathbb C\setminus D(r_0)$, and  $f$  vanishes on a sequence  ${\tt Z}=\{{\tt z}_k\}_{k=1,2,\dots}\subset \mathbb C\setminus D(r_0)$.
	If $v\in \text{\rm sbh}^+\bigr(\mathbb C \setminus D(r_0)\bigr)$, $	\lim\limits_{z\to \infty} v(z){=}0$,
and 
\begin{equation*}
	\int_{\mathbb C\setminus D(r_0)}  v \,d {\nu}_M<+\infty, \quad\text{but}\quad
		\sum_{k}v({\tt z}_k)=+\infty,
\end{equation*}
then  $f=0$.}
\medskip 

\noindent
\textsc{Corollary 2.} {\it   Let $r_0>0$, and   $v\in \text{\rm sbh}^+\bigr(\mathbb C \setminus D(r_0)\bigr)$
   with the Riesz measure $\nu_v$ satisfies the condition  $	\lim\limits_{z\to \infty} v(z){=}0$. 
If an entire  function $f\neq 0$ vanishes on a sequence  ${\tt Z}=\{{\tt z}_k\}_{k=1,2,\dots }\subset \mathbb C$,
and $f$ satisfies the condition 	$ 	\int\limits_{\mathbb C \setminus D(R_0)} \log |f| \, d \nu_v<+\infty $, then 
$$
\sum\limits_{{\tt z}_k \in D\setminus D_0} v({\tt z}_k)<+\infty.
$$}
\medskip

\noindent
\textsc{{\bf 4.} Radial versions.}  Let $r_0\in (0,+\infty)$. Let $M\colon (r_0,+\infty) \to \mathbb R$ be an increasing continuous function. 
Further assume that there exists  its left-hand derivative $M'_{\text{\tiny\rm left}}$ on $(r_0, +\infty)$
and the function $t\mapsto tM'_{\text{\tiny\rm left}}(t)$ is increasing on $(r_0,+\infty)$.

Let  $q\colon [r_0,+\infty)\to \mathbb R$ be a bounded positive decreasing function. Further assume that 
$\int\limits_{r_0}^{+\infty} {q(t)}\, \dfrac{d t}{t}<+\infty .	$

Let $f\neq 0$ be an entire function such that $f$  vanishes on a sequence  ${\tt Z}=\{{\tt z}_k\}_{k=1,2,\dots }\subset \mathbb C\setminus \overline{D(r_0)}$.
\medskip

\noindent
\textsc{Corollary 3.} {\it    If\/  $\bigl|f(z)\bigr|\leq \exp M\bigl(|z|\bigr)$ for all\/ $|z|>r_0$ and
\begin{equation*}
	\int_{r_0}^{+\infty} q(t)M'_{\text{\tiny\rm left}}(t) \, dt<+\infty,
\end{equation*}
then
\begin{equation*}
	\sum\limits_{|{\tt z}_k |>r_0} \int_{|{\tt z}_k|}^{+\infty} \frac{q(t)}{t} \, dt <+\infty.
\tag{\rm q}
\end{equation*}}

\noindent
\textsc{Corollary 4.} {\it  If 	$\int\limits_{r_0}^{+\infty}\Bigl(\frac{1}{2\pi} \int\limits_0^{2\pi}\log \bigl|f(te^{i\theta})\bigr| \, d\theta\Bigr) \, d q(t)  <+\infty$, then the relation\/ {\rm (q)} is fulfilled.}
\medskip

These Corollaries 3,4 improve a part of results from [2]--[8], [12]. Explicit non-radial cases are much more complicated and will be discussed in detail in another place.

Supported by RFBR--РФФИ (grant no. 13-01-00030a).

\medskip
\centerline{\bf References}
\smallskip

1.\; Ransford\,Th. {\it Potential Theory in the Complex Plane}~-- Camb\-r\-i\-d\-ge:  Cambridge University Press, 1995.

2.\;Хабибуллин\,Б.\,Н. {\it Полнота систем экспонент и множества единс\-т\-в\-е\-н\-н\-о\-с\-ти\/} (издание $4^{\text{\underline{ое}}}$, дополненное). 
-- Уфа: РИЦ БашГУ. 2012.\\
{http://www.researchgate.net/profile/Bulat{\_}Khabibullin/contributions}

3.\;Khabibullin\,B.\,N. {\it Sets of uniqueness in spaces of entire functions of a sin\-gle variable} //~Mathematics of the USSR--Izvestiya.~-- 1992. 
V.~39.~-- no.~2.~-- P.~1063--1084. --- Изв. АН СССР.  Сер. матем.~-- 1991.~-- Т. 55.~-- №~5.~-- C.~1101--1123.

4.\;Khabibullin\,B.\,N. {\it Distribution of zeros of entire functions and the balayage.\/} 
Diss. of Dr. of ph.-math. Sci.~-- 1993. Kharkiv. 

5.\;Khabibullin\,B.\,N. {\it Dual representation of superlinear fu\-n\-c\-tionals and its applications in function theory.} II // 
Izvest\-i\-ya: Math\-e\-m\-a\-t\-ics.~-- 2001.~-- V.~65.~-- no.~5.~-- P.~1017--1039. --- Изв. РАН. Сер. ма\-т\-ем.~-- 2001.~--
 Т. 65.~--  № 5.~-- С. 167-190.

6.\;Khabibullin\,B.\,N. {\it Zero sequences of holomorphic functions, representation of meromorph\-ic functions, and harmonic minorants} //
Sbornik: Mathematics.~--  2007.~-- 198\,:\,2.~-- P.~261--298. --- Матем. сб.~-- 2007.~--  Т.~198.~-- №~2.~-- С.~121--160. 

7.\;Khabibullin\,B.\,N.,  Khabibullin\,F.\,B., Cherednikova\,L.\,Yu. {\it Zero subsequences for classes of holomorphic functions: stability and the entropy of arcwise connectedness.\/} II // St. Petersburg Ma\-th. Journal.~-- 2009.~--  20\,:\,1.~-- P.~131--162. --- Алгебра и анализ. 2008. Т.~20. №~1. С.~190--236. 

8.\;{Khabibullin\,B.\,N.} {\it Sequences of non-uniqueness for weight spaces of holomorphic functions} // Russian Mathematics (Izvestiya VUZ. Matematika)~--
 2015~-- 59\,:\,4.~-- P.~63--70. --- Известия вузов. Математика.~-- 2015.~-- № 4.~-- С.~75--84.

9.\;Matsaev\,V., Sodin\,M. {\it Distribution of Hilbert transforms of measures} //  Geom. funct. anal.~-- 2000.~--
V.~10.~-- P.~160--184. 

10.\;Khabibullin\,B.\,N., Talipova\,G.\,R.,  Khabibullin\,F.\,B. {\it Zero su\-b\-s\-e\-quences for Bernstein's spaces and the completeness of ex\-p\-o\-n\-e\-ntial systems in spaces of functions on an interval} // St. Petersburg Mathematical Journal.~--  2015~-- 26\,:\,2.~-- P.~319--340. ---
 Алгебра и анализ.~-- 2014.~-- Т.~26, №~2.~--  С.~185--215. 

11.\;Talipova\,G.\,R.,  Khabibullin\,B.\,N. {\it Sequences of uniqueness for classes of entire functions of exponential type with restrictions on the real axis} //
Bulletin of  Bashkir University.~-- 2015.~-- V.~20.~-- no.~1~-- P.~1--5 --- Вестник Башкирского ун-та.~-- 2015.~--
 Т. 20.~-- №~1.~-- С.~1--5.

12.\;Bykov\,S.\,V., Shamoyan\,F.\,A.  {\it On zero of entire functions with majorant of infinite order} //
St. Petersburg Mathematical Journal.~--  2010~-- 21\,:\,6.~--  P.~893--901.  --- Алгебра и анализ --- 2009. --- Т.~21.~-- №~6.~-- С.~66--79. 

\end{document}